\RequirePackage{ifpdf}
\ifpdf 
\documentclass[pdftex]{sigma}
\else
\documentclass{sigma}
\fi

\usepackage{xypic}

\numberwithin{equation}{section}
\numberwithin{theorem}{section}

\newtheorem{Proposition}[theorem]{Proposition}
\newtheorem{Question}[theorem]{Question}
\newtheorem{Conjecture}[theorem]{Conjecture}
\newtheorem{Corollary}[theorem]{Corollary}
{\theoremstyle{definition} \newtheorem{Definition}[theorem]{Definition}
}

\makeatletter
\def\revddots{\mathinner{\mkern1mu\raise\p@
    \vbox{\kern7\p@\hbox{.}}\mkern2mu
    \raise4\p@\hbox{.}\mkern2mu\raise7\p@\hbox{.}\mkern1mu}}
\makeatother

\newcommand{\blambda}{{\boldsymbol{\lambda}}}
\newcommand{\bnu}{{\boldsymbol{\nu}}}
\newcommand{\bmu}{{\boldsymbol{\mu}}}

\newcommand{\bA}{\boldsymbol{A}}

\newcommand{\Cset}{\mathbb{C}}
\newcommand{\hC}{\hat{\Cset}}

\newcommand{\Stab}{\operatorname{Stab}}
\newcommand{\M}{\operatorname{M}}

\newcommand{\I}{\operatorname{I}}
\newcommand{\W}{\operatorname{W}}
\newcommand{\V}{\operatorname{V}}
\newcommand{\R}{\operatorname{R}}
\renewcommand{\P}{\operatorname{P}}

\newcommand{\B}{\operatorname{B}}

\newcommand{\tnu}{\tilde{\nu}}

\newcommand{\hU}{{\hat{U}}}
\newcommand{\hp}{{\hat{p}}}
\newcommand{\hlambda}{{\hat{\lambda}}}
\newcommand{\hz}{\hat{z}}

\newcommand{\hnu}{\hat{\nu}}
\newcommand{\hmu}{\hat{\mu}}
\newcommand{\hblambda}{\hat{\blambda}}

\newcommand{\hR}{\hat{R}}

\newcommand{\lspan}{\operatorname{span}}
\newcommand{\ord}{\operatorname{ord}}
\newcommand{\SL}{\operatorname{SL}}
\newcommand{\PSL}{\operatorname{PSL}}

\newcommand{\cP}{\mathcal{P}}
\newcommand{\cG}{\mathcal{G}}

\newcommand{\pmat}[1]{\begin{pmatrix}#1\end{pmatrix}}

\newcommand{\supth}{^{\mathrm{th}}}

\begin{document}


\renewcommand{\PaperNumber}{107}

\FirstPageHeading

\ShortArticleName{On Projective Equivalence of Univariate Polynomial Subspaces}

\ArticleName{On Projective Equivalence\\ of Univariate Polynomial Subspaces}

\Author{Peter CROOKS~$^\dag$ and Robert MILSON~$^\ddag$}

\AuthorNameForHeading{P.~Crooks and R.~Milson}

\Address{$^\dag$~Department of Mathematics, University of Toronto,
 Toronto, Ontario, Canada M5S 2E4}
\EmailD{\href{mailto:peter.crooks@utoronto.ca}{peter.crooks@utoronto.ca}}

\Address{$^\ddag$~Department of Mathematics and Statistics, Dalhousie
  University,\\
  \hphantom{$^\ddag$}~Halifax, Nova Scotia, Canada B3H 3J5}
\EmailD{\href{mailto:rmilson@dal.ca}{rmilson@dal.ca}}

\ArticleDates{Received June 05, 2009, in f\/inal form December 03, 2009;  Published online December 06, 2009}

\Abstract{We pose and solve the equivalence problem for subspaces of $\cP_n$,
  the $(n+1)$ dimensional vector space of univariate polynomials of
  degree $\leq n$.  The group of interest is $\SL_2$ acting by
  projective transformations on the Grassmannian variety $\cG_k\cP_n$
  of $k$-dimensional subspaces.  We establish the equivariance of the
  Wronski map and use this map to reduce the subspace equivalence
  problem to the equivalence problem for binary forms.}


\Classification{14M15; 15A72; 34A30; 58K05}

\section{Introduction}

Given a representation of a Lie group $G$ on a f\/inite-dimensional
vector space $V$, it is natural to consider the associated action of
$G$ on $\cG_k V$, the Grassmannian variety of $k$-dimensional
subspaces of $V$.  This non-linear $G$ action is realized via the
so-called Pl\"ucker embedding of $\cG_k V$ into the projectivization
of $\Lambda^k V$, the $k\supth$ exterior product of $V$.  The
following question is of central importance:
\begin{quote}
  \textit{When are two subspaces $U_1,U_2\in \cG_kV$ equivalent up to
    a $G$-action?}
\end{quote}

The present article is devoted to a specif\/ic instance of this general
equivalence problem. Here, the group of interest is $\SL_2$ acting on
\begin{equation*}
  \cP_n(z): =\lspan \{ 1,z,z^2,\ldots, z^n \},
\end{equation*}
by fractional linear transformations.  Classical invariant
theory~\cite{hilbert,olver99} is concerned with the equivalence problem
for one-dimensional subspaces; i.e., binary forms, def\/ined up to scalar
multiple.  We are interested in the following generalization.
\begin{Definition}
  \label{def:projequiv}
  Let $U,\hU\in \cG_k\cP_n$ be polynomial subspaces with respective
  bases $p_i, \hp_i\in \cP_n$, $i=1,\ldots, k$.  We say that the two
  subspaces are projectively equivalent, and write $U\sim \hU$ if there
  exist a fractional linear transformation
  \begin{equation*}
    z = \frac{a\hz + b}{c\hz+d},\qquad ad-bc=1,
  \end{equation*}
  and an invertible $k\times k$ matrix $A_{ij}$ such that
  \[ \hp_i(\hz)  = (c\hz+d)^n \sum_{j=1}^k A_{ij}
  p_j \left(\frac{a\hz+b}{c\hz+d}\right ),\qquad
  i=1,\ldots,k.\]
\end{Definition}

\looseness=1
As we shall see, the equivalence problem for polynomial subspaces
furnishes interesting analogues and extensions of certain concepts
from the classical invariant theory of binary forms.

\looseness=1
Our original motivation for studying the equivalence problem for
polynomial subspaces relates to the theory of dif\/ferential operators
with invariant polynomial subspaces~\cite{kamfink}.  The two
fundamental questions in this area are: (1)~to characterize polynomial
subspaces that are left invariant by a 2nd order dif\/ferential
operator; and (2)~to characterize inf\/inite polynomial f\/lags preserved
by a 2nd order operator.  The f\/irst question is central to the theory
of quasi-exact solvability in quantum mechanics~\cite{gkm1}.  The
second question has interesting applications to orthogonal polynomials~\cite{gkm3}.  Invariance with respect to a dif\/ferential operator is an
$\SL_2$-invariant property, and further progress depends crucially on
advancing the understanding of invariant properties of polynomial
subspaces.  The present paper gathers the necessary exposition and
theorems. We plan to return to the question of dif\/ferential operators
in subsequent publications.

The main tool in our analysis of the subspace equivalence problem is
the Wronskian operator.  Let $p_1,\ldots, p_k\in \cP_n$ be linearly
independent polynomials and $W = \W(p_1,\ldots, p_k)$ the
correspon\-ding Wronskian determinant. It isn't hard to show that $\deg
W \leq k\ell$, where $\ell=n+1-k$ is the codimension of a
$k$-dimensional subspace of $\cP_n$.  Since equality holds
generically, we can def\/ine the linear Wronski map $\W:\Lambda^k \cP_n
\to \cP_{k\ell}$, and by restriction its non-linear
analogue\footnote{We use the same symbol $\W$ for both maps and for
  the Wronskian operation.  Ambiguity isn't a danger, because the
  particular choice of usage is indicated by the context.}  $\W:\cG_k
\cP_n \to \cG_1 \cP_{k\ell}$. For details, see Proposition~\ref{prop:degW} and the extensive discussion in~\cite{scherbak}. The Wronski map is intimately related to the Schubert
calculus~\cite{fulton} and the enumeration of rational curves in
projective space~\cite{khsott03}, and has received a great deal of
recent attention owing to its connection with the exactly solvable
Gaudin model and the Shapiro conjecture; see~\cite{eremenko, scherbak,  mukhin} and the references therein.

Key to our analysis is the fact that the the Wronski map is $\SL_2$
equivariant.  From the representation theory point of view,
$\cP_{k\ell}$ is just the highest weight component of the
decomposition of $\Lambda^k\cP_n$ into irreducible $\SL_2$
representations.  Reference~\cite{chipalkatti} explores this point of
view, as well as the relationship of the Wronskian to classical
invariant theory.  In classical terminology, the Wronskian is just a
particular combinant, that is a joint covariant of polynomials $p_1,
\ldots, p_k$ that enjoys a determinantal transformation law with
respect to a change of basis.  The projections of~$\Lambda^k\cP_n$
onto the irreducible $\SL_2$ components is accomplished by a sequence
of Wronskian-like combinants $C_i: \Lambda^k \cP_n \to \cP_{k\ell-2i}$
where $i=0,2,\ldots, k$ and where $C_0$ is just the ordinary
Wronskian.  It is shown in \cite{chipalkatti} that the mapping
\[ \lspan(p_1,\ldots, p_k) \mapsto [(C_0,C_2,\ldots, C_k)] \] def\/ines
an equivariant projective embedding of the Grassmannian $\cG_k\cP_n$
into a projective space of suf\/f\/iciently high dimension.  Therefore, in
principle, the equivalence problem for polynomial subspaces is fully
solved by considering the joint $\SL_2$ equivalence problem~\cite{olverji01} for the
combinants $(C_0, C_2,\ldots, C_k)$.

Our main observation is that the equivalence problem for polynomial
subspaces can be solved by considering just $C_0$, the ordinary
Wronskian.  The equivariance of the Wronski map means that if two
subspaces are projectively equivalent, then so are their respective
Wronskians.  The converse, of course, is not true.  The consequence of
discarding the other combinants is that the resulting mapping is no
longer one-to-one. Indeed, the degree of the Wronski map $U\mapsto
W(U) = C_0$ is of central concept in the Schubert calculus and has an
important combinatorial interpretation~\cite{fulton,griffiths}.
However, we can assert the following.
\begin{Proposition}
  \label{prop:genequiv}
  Let $U, \hU\in \cG_k\cP_n$ be polynomial subspaces such that $W(\hU) =
  g\cdot W(U)$ where $g\in \SL_2\Cset$ is a projective
  transformation.  Let us furthermore suppose that $\Stab W(U)$, the
  subgroup of projective transformation that leaves $W(U)$ invariant,
  is trivial.  Then $U\sim \hU$ if and only if $\hU = g\cdot U$.
\end{Proposition}

Thus, in the generic case, we can decide the equivalence problem for
polynomial subspa\-ces~$U$,~$\hU$ as follows.  First we decide whether the
corresponding Wronskians are equivalent.  If the Wronskians are
inequivalent, then so are $U$ and $\hU$. However, if the Wronskians
are equivalent, \emph{and if we can determine the fractional linear
  transformation that relates~$W(U)$ and~$W(\hU)$}, then we can also
decide the question of equivalence of the subspaces $U$ and $\hU$.

\looseness=-1
The classical approach to the equivalence of binary forms is based on
an integrity basis~-- a~f\/inite generating set of invariants or
covariants~\cite{kungrota}.  Here, we follow a dif\/ferent, more modern
approach based on moving frames introduced in~\cite{olver99}.  It
turns out that if a polynomial has a~f\/inite symmetry group, then it
can be fully characterized by two absolute rational covariants.
In~\cite{BO00} the equivalence problem for binary forms is
reduced to a procedure that tests the compatibility of these two
covariants.  In the case of a positive outcome of the equivalence
test, the procedure actually determines the projective transformation
$g\in \SL_2\Cset$ that relates the binary forms in question.  If we
apply the algorithm of~\cite{BO00} to Wronskians $W(U)$ and $W(\hU)$
and f\/ind a~transformation $g$ such that $W(\hU)= g\cdot W(U)$, then,
in the generic case, we can test for the equivalence of~$U$ and~$\hU$
as per the above Proposition.  The case where the Wronskian has an
inf\/inite symmetry group can also be handled by means of covariants.
See Section~\ref{sect:eqprob} for the details.

Here is the outline of the rest of the paper. Apolarity is another key
notion in the classical invariant theory~\cite{kungrota}.  The
corresponding, and closely related, notion for polynomial subspaces is
equivariant duality between dimension and codimension.  The duality is
def\/ined and discussed in Section~\ref{sect:dual}.  Invariance of the
Wronski map with respect to the apolar duality is established in
Theorem~\ref{thm:Wdual}.  This theorem is far from surprising, and the
proof is straightforward; however, we were unable to locate the
theorem in the literature.

Section~\ref{sect:partition} gathers the necessary background on
partitions, reduced row echelon form, and canonical reduced bases.
Section \ref{sect:inv} is devoted to the proof of the equivariance of
the Wronski map.  These sections are largely expository, but we
include them for the convenience of the reader.  An explicit proof of
equivariance of the Wronskian can be found in~\cite{chipalkatti}.  In the present paper, we present two alternative,
elementary proofs of equivariance.  The f\/irst proof is novel (again, to our
best knowledge,) while the second builds on a proof sketch in~\cite{eremenko}.

Section~\ref{sect:eqprob} is devoted to the equivalence problem for
polynomial subspaces.  There we present a~number of novel (to our best
knowledge) results.  Theorem \ref{thm:genequiv} concerns the reduction
of the the equivalence problem for subspaces to the equivalence
problem for binary forms in the case of a~f\/inite symmetries.
Theorem \ref{thm:onetworoots} deals with the non-generic case of
inf\/inite symmetry groups. We prove that a subspace~$U$ has an
inf\/inite symmetry group if and only if, up to a~fractional linear
transformation, it is spanned by monomials.  Equivalently, we can
characterize such subspaces by a condition involving the Hessian of
the Wronskian and another covariant.

Theorem \ref{thm:npower} concerns the question of top powers.  The
classical Waring problem for polynomials asks for the smallest integer
$s\geq 1$ such that a $p\in \cP_n$ can be written as $s$-term sum of
$n$th powers.  An analogous problem for polynomial subspaces is
formulated in~\cite{chipger}.  The authors of that paper also
study the following, related question: for $s\geq k$ describe the
closure of the locus in $\cG_k(\cP_n)$ formed by subspaces that are
spanned by at most $s$ top powers.  We limit ourselves to the case
$s=k$ and focus on the following question: to describe a decision
procedure that determines whether a $k$-dimensional polynomial
subspace $U\in \cG_k\cP_n$ admits a basis of $n$th, or top powers.  An
algorithmic criterion for this condition is presented in Theorem~\ref{thm:npower}.  Finally, in Proposition~\ref{prop:lthpower} we
present a striking necessary condition for the existence of a basis of
top powers and conjecture that this condition is also suf\/f\/icient.  A
proof of the conjecture is given for the cases of codimension 1 and
codimension 2 subspaces.

\section{Apolar duality}
\label{sect:dual}
For convenience, we work over $\Cset$, the f\/ield of complex numbers.
The group $\SL_2\Cset$ is the group of unimodular transformations of
$\Cset^2$:
\begin{equation*}
  \pmat{\hat{x}\\\hat{y}} = \pmat{a&b\\c&d} \pmat{x\\y},\qquad ad-bc=1.
\end{equation*}
Setting $z=x/y$ we obtain the fractional linear
transformation
\begin{equation*}
  z = \frac{a\hz+b}{c\hz+d},\qquad ad-bc=1
\end{equation*}
of complex projective space $\Cset\P^1$. We identify the latter with
the extended complex plane, $\hC := \Cset \cup \{ \infty\}\cong
\Cset\P^1$.  The set of all such transformations forms the quotient
group $\PSL_2\Cset \simeq \SL_2\Cset/\{\pm I\}$.  Fractional linear
transformations also realize the $(n+1)$ dimensional irreducible
representation of $\SL_2\Cset$ by means of the following action on
$\cP_n$:
\begin{equation*}
  p(z)\mapsto \hp(\hz)  = (c \hz+d)^n\, p(z),\qquad p \in   \cP_n.
\end{equation*}

Let $\cG_{k}\cP_n$ denote the Grassmannian variety whose elements are
$k$-dimensional subspaces of the $n+1$ dimensional vector space
$\cP_n$.
In our analysis of this problem, we introduce the non-degenerate
bilinear form $\gamma:\cP_n\times \cP_n\to \Cset$, def\/ined by the
relations:
\begin{equation}
  \label{eq:gammabracket}
  n!\, \gamma\left(\frac{z^j}{j!},\frac{z^k}{k!}\right)= \begin{cases}
  \displaystyle (-1)^j, & \text{if} \ \ j+k=n,\\
  0, & \text{otherwise.}
\end{cases}
\end{equation}
Equivalently, we can write
\begin{equation*}
  \gamma^{-1}=\sum_{j=0}^n (-1)^{j}
  \binom{n}{j} z^j \otimes z^{n-j}.
\end{equation*}
Note that $\gamma$ is symmetric if $n$ is even, and skew-symmetric if
$n$ is odd.  Also note that, up to sign, $\gamma(p,q)$ is the
classical apolar invariant of two polynomials $p,q\in \cP_n$, also
known as the lineo-linear invariant~\cite[Section~5]{kungrota}.  This
remark establishes the following result.  A
direct proof is also given in~\cite{gkm2}.
\begin{Proposition}
  The above-defined bilinear form $\gamma$ is $\SL_2\Cset$-invariant.
\end{Proposition}

\begin{Definition}
  For a polynomial subspace $U\in \cG_k\cP_n$ we def\/ine
  \begin{equation*}
    U^* = \{ u\in \cP_n \colon \gamma(u,v)=0\text{ for all } v\in U\},
  \end{equation*}
  and refer to $U^*$ as the apolar dual.
\end{Definition}
\begin{Proposition}
  The apolar duality mapping $\cG_k \cP_n\to\cG_{n+1-k}\cP_n$ given by
  \[U\mapsto U^*,\qquad U\in \cG_k\cP_n\] is $\SL_2$-equivariant.
\end{Proposition}

Also observe that if $p_1,\ldots, p_k$ is a basis of $U$, then a basis
for $U^*$ is obtained by solving a~system of $k$ homogeneous linear
equations in $n+1$ variables, namely{\samepage
\begin{equation*}
  \gamma(p_j,q) = 0,\qquad j=1,\ldots, k
\end{equation*}
where $q\in \cP_n$ is the unknown.}

Classically, two polynomials of degree $n$ that satisfy
$\gamma(p,q)=0$ are called apolar\footnote{The classical def\/inition of
  apolarity for polynomials of dif\/ferent degrees $m>n$ involves a
  bilinear bracket $\cP_m \times\cP_n \mapsto \cP_{m-n}$.  The
  classical bracket matches our form $\gamma$ only in the case where
  $m=n$.}.  The vector space of binary forms apolar to a given $p\in
\cP_n$ is a fundamental object in classical invariant theory and is
closely related to the question of sums of $n$th powers; see Kung and
Rota~\cite{kungrota} for details.  Indeed we have the following
classical result.

\begin{Proposition}
  Suppose that $p\in \cP_n$ has $n$ distinct, finite roots.  Then, the
  $n$ dimensional dual vector space $\lspan\{ p\}^*$ is spanned by $n$th
  powers of the form $(z-r_i)^n$ where $r_1,\ldots, r_n$ are the roots
  of $p(z)$.
\end{Proposition}

We will return to the question of $n$th powers in the f\/inal section of
our paper.

\section{Partitions, diagrams, and reduced echelon form}
\label{sect:partition}

\subsection{Partitions}
\label{sub:partitions}

A partition of an integer $N\geq 0$ is a f\/inite multi-set\footnote{The
  empty set corresponds to $N=0$.} of positive integers that sum to
$N$.  For our purposes, it will be convenient to def\/ine a partition as
a non-increasing sequence $\blambda=\{\lambda_i\}$ of non-negative
integers $\lambda_1\geq \lambda_2\geq \cdots \geq 0$ such that
\begin{equation*}
  N(\blambda) = \sum_i\lambda_i <\infty.
\end{equation*}
We will refer to  $N$ as the partition \emph{total}, and call
\begin{gather*}
  h(\blambda) = \max \{ i: \lambda_i>0 \},\qquad
  w(\blambda)=\lambda_1,
\end{gather*}
respectively, the height and the width.  The former names the
cardinality of the multi-set, while the latter names the largest
element of the multi-set. The height and width terminology derives
from the fact that a partition can be visualized by means of a
Ferrer's diagram: $\lambda_1$ dots in the f\/irst row, $\lambda_2$ dots
in the 2nd row, etc.  The trivial partition is one where
$\lambda_i=0$ for all $i$; this corresponds to the empty multi-set.
By agreement, the trivial partition has $h=w=N=0$.

\begin{Definition}
  \label{def:conjpartchar}
  Let $\blambda=\{\lambda_i\}$ be a partition.  For $j\geq 1$, def\/ine
  $\lambda^*_j$ to be the cardinality of the set $\{ i\geq 1:
  \lambda_i \geq j\}$.  Call the resulting sequence
  $\blambda^*=\{\lambda^*_j\}$ the conjugate partition of $\blambda$.
\end{Definition}
\begin{Proposition}
  \label{prop:conjpartchar}
  Let $\blambda$ be a partition of height $h$ and width $w$.  The
  conjugate partition $\blambda^*=\{\lambda^*_j\}$ is the unique
  partition of height $w$ and width $h$ such that for all $i,j\geq 1$
  we have
  \begin{equation*}
    j\leq \lambda_i \quad\text{if and only if}\quad i\leq \lambda^*_j.
  \end{equation*}
  Furthermore, conjugate partitions have the same total: $N(\blambda)=
  N(\blambda^*)$.
\end{Proposition}

In other words, the Ferrer's diagram of the conjugate partition
$\blambda^*$ is obtained by transposing the left-aligned Ferrer's
diagram of $\blambda$; columns become rows, and vice versa.  In the
next section, we will characterize the apolar duality relation between
polynomial subspaces using the conjugate partition.  First, we
characterize conjugate partitions in terms of complementary subsets.

From now on, we f\/ix integers $k,\ell\geq 0$  and set
\begin{gather*}
  n=k+\ell-1.
\end{gather*}
We will say that a partition is \emph{$k\times \ell$ bounded} if
$h\leq k$ and $w\leq \ell$.  In other words, for $k,\ell \geq 1$, the
Ferrer's diagram of a $k\times \ell$ bounded partition is a subset of
the discrete $k\times \ell$ rectangle.  Let~$\B_{k,\ell}$ denote the
set of all $k\times \ell$ bounded partitions.  Thus, the conjugation
operation $\blambda\to \blambda^*$ def\/ines a bijection of
$\B_{k,\ell}$ and $\B_{\ell,k}$.  We say that a partition
$\blambda=\{\lambda_i\}$ is \emph{rectangular} if every non-zero~$\lambda_i$ is equal to $\lambda_1$. The $k\times \ell$ rectangular
partition is the ``largest'' element of $\B_{k,\ell}$.

Let $\P_{k,\ell}$ denote the set of all strictly increasing
$k$-element sequences\footnote{When $n=-1$ we are
  speaking of subsets of the empty set.  When $n=0$ we are referring
  to subsets of $\{0\}$.} in $\{0,1,\ldots, n\}$.  Thus, $\bnu\in
\P_{k,\ell}$ refers to a $k$-element sequence $0\leq \nu_1<\cdots <
\nu_k\leq n$.  We now describe two useful bijections of $\B_{k,\ell}$
and $\P_{k,\ell}$.  For a bounded partition $\blambda\in \B_{k,\ell}$
def\/ine $\blambda^+,\blambda^-\in \P_{k,\ell}$ by setting
\begin{gather}
  \label{eq:lambda+def}
  \lambda^+_i  = \ell-1-\lambda_i+i,\qquad i=1,\ldots, k,\\
   \lambda^-_i  = \lambda_{k+1-i}+i-1
  = n-\lambda^+_{k+1-i},\qquad i=1,\ldots, k. \label{eq:lambda-def}
\end{gather}
Observe that $\lambda^+_1 = \ell-\lambda_1$, $\lambda^+_k =
n-\lambda_k$.  Similarly, $\lambda^-_1 = \lambda_k$, $\lambda^-_k =
\lambda_1+k-1$.  Therefore the mappings $\blambda\mapsto \blambda^+$
and $\blambda\mapsto \blambda^-$ are bijective.

In addition to conjugation, bounded partitions enjoy another kind of
duality, one that relates~$\blambda^+$ and~$\blambda^-$.  Given a
bounded partition $\blambda\in \B_{k,\ell}$, set
\begin{gather*}
  \lambda'_i  := \ell - \lambda_{k+1-i},\qquad i=1,\ldots, k,\qquad
  \lambda'_i  := 0,\qquad i> k
\end{gather*}
and call the resulting partition $\blambda'=\{\lambda'_i\}$ the
$k\times l$ complement of $\blambda$.  Equivalently, the Ferrer's
diagram of $\blambda'$ is the reversed complement in a $k\times \ell$
rectangle of the Ferrer's diagram of $\blambda$.  Also, observe that
\begin{equation*}
  N(\blambda') = k\ell - N(\blambda).
\end{equation*}

\begin{Proposition}
  \label{prop:partcomp}
  Bounded partitions $\blambda,\bmu\in \B_{k,\ell}$ are complementary
  partitions, $\bmu = \blambda'$, if and only if $\blambda^+=\bmu^-$.
\end{Proposition}

We now come to the main result of this section.
\begin{theorem}
  \label{thm:conjpart}
  Bounded partitions $\blambda\in \B_{k,\ell}$ and $\hat{\blambda}\in
  \B_{\ell,k}$ are conjugate partitions if and only if $\blambda^+$
  and $\hat{\blambda}{}^{\!-}$ are complementary subsets of
  $\{0,1,\ldots, n\}$.
\end{theorem}

\subsection{Row reduced echelon form}

Let $\blambda\in \B_{k,\ell}$ be a bounded partition.  We will say
that a $k\times \ell$ matrix $\{ a_{ij} \}$ where $1\leq i\leq k$ and
$1\leq j \leq \ell$ is
$\blambda$-bounded if
$a_{ij} = 0$ for all $j\leq \ell-\lambda_i$.

Next we describe a useful bijection between the set of
$\blambda$-bounded matrices and matrices in reduced row echelon form
(RREF).  Let $\R_{k,\ell}$ denote the set of $k\times (k+\ell)$
matrices that have rank~$k$ and that are in RREF.  When dealing with
matrices in $\R_{k,\ell}$ we adopt the convention that $n=k+\ell-1$,
while $1\leq i\leq k$ and $0\leq j \leq n$ serve as the row and column
indices, respectively.  These index conventions are convenient because
of subsequent applications to polynomial subspaces; the column indices
$\{0,1,\ldots, n\}$ enumerate the possible degrees of a polynomial in
$\cP_n$.

Every matrix in $\R_{k,\ell}$
can be uniquely block-partitioned as follows
\begin{equation}
  \label{eq:RREF}
  \begin{pmatrix}
    \I_1 & B_{11} & 0 & B_{12} & \cdots & 0& B_{1r} \\
    0   & 0  & \I_2   & B_{22} & \cdots & 0& B_{2r} \\
 \vdots & \vdots & \vdots & \vdots & \ddots & \vdots & \vdots \\
 0 & 0 & 0 & 0 & \cdots & \I_r & B_{rr}
  \end{pmatrix},
\end{equation}
where the submatrices $B_{pq}$, $1\leq p\leq q\leq r$ are blocks
of size $k_p \times \ell_q$ and
the blocks $\I_p$ denote a~$k_p\times k_p$ identity matrix. Here
 \begin{gather*}
   k_1+\cdots + k_r = k,\qquad
    \ell_1 + \cdots + \ell_r = \ell,
\end{gather*}
with $k_1,\ell_r\geq 0$  and  all other
$k_i$, $\ell_j$ strictly positive.
The possibilities $k_1=0$ and $\ell_r=0$ correspond to 3 degenerate
subcases, shown below:
\begin{gather}
  \label{eq:degen1}
      \begin{pmatrix}
     0  & \I_2   & B_{22} & \cdots & 0& B_{2r} \\
     \vdots & \vdots & \vdots & \ddots & \vdots & \vdots \\
     0 & 0 & 0 & \cdots & \I_r & B_{rr}
  \end{pmatrix},\qquad k_1=0,\quad \ell_r >0;\\
  \label{eq:degen2}
   \begin{pmatrix}
    \I_1 & B_{11} & 0 & B_{12} & \cdots & 0 \\
    0   & 0  & \I_2   & B_{22} & \cdots & 0 \\
    \vdots & \vdots & \vdots & \vdots & \ddots & \vdots \\
    0 & 0 & 0 & 0 & \cdots & \I_r
  \end{pmatrix},\qquad k_1 > 0,\quad \ell_r=0;\\
  \label{eq:degen3}
   \begin{pmatrix}
    0  & \I_2   & B_{22} & \cdots & 0 \\
    \vdots & \vdots & \vdots & \ddots & \vdots \\
     0 & 0 & 0 & \cdots &  \I_r
  \end{pmatrix},\qquad k_1= \ell_r=0.
\end{gather}
Also note that $\R_{k,0}$ is a singleton set consisting of the $k\times
k$ identity matrix.

Deleting the pivot columns from a matrix in RREF results in a
$\blambda$-bounded matrix for a~certain bounded partition
$\blambda\in \B_{k,\ell}$.
\begin{Proposition}
  \label{prop:lpartRREF}
  Let $\blambda\in \B_{k,\ell}$ be a bounded partition. Let
  $\bnu=\blambda^+$ and let $\bmu$ be the increasing enumeration of
  the complement $\{0,1,\ldots, n\}\backslash\{\nu_i\}$. Let $\{a_{ij}\}$ be a
  $\blambda$-bounded matrix, and let $\bA=(A_{ij})$ be the $k\times (k+\ell)$
  matrix defined by
  \begin{gather*}
     A_{i\mu_j} = a_{ij},\qquad i=1,\ldots, k,\quad j=1,\ldots, \lambda_i,\\
     A_{i\nu_i} = 1,
  \end{gather*}
  with all other entries zero.  Then, $A$ is in RREF; $\{\nu_i\}$ is
  the enumeration of the pivot columns of $A$; and $\{\mu_j\}$ is the
  enumeration of the non-pivot columns.  Conversely, every $A\in
  \R_{k,\ell}$ admits a unique such representation in terms of a
  $\blambda$-bounded  matrix for some $\blambda\in \B_{k,\ell}$.
\end{Proposition}

 The above construction associates the trivial partition to
the most extreme form of \eqref{eq:degen3}, the matrix composed of
$\ell$ zero columns followed by the $k\times k$ identity matrix.

\subsection{The shape of a polynomial subspace}

Next, we apply the above results on reduced row echelon form to
polynomial subspaces.  Fix integers $n\geq 0$ and $1\leq k\leq n+1$
and set $\ell =n+1-k$.  For a non-zero polynomial $p\in \cP_n$ and
$b\in \Cset$ we def\/ine
\begin{equation*}
  \ord_bp = \min\{ j\geq 0 : p^{(j)}(b)\neq 0 \}.
\end{equation*}
As per the usual convention, $\ord_b 0 = \infty$.
For a polynomial subspace $U$, we def\/ine
\begin{equation*}
  \ord_b U = \min \{ \ord_b p : p \in U \}.
\end{equation*}
Equivalently, if $U$ is non-trivial, $\nu=\ord_b U$ is the largest
integer such that $(z-b)^\nu$ divides \emph{all} $p(z)\in U$.

\begin{Proposition}
  \label{prop:ordU}
  Let $U$ be a polynomial subspace and let $b\in \Cset$.  Then, the
  subspace
  \begin{equation*}
    U'=\{ p\in U: \ord_b p > \ord_b U \}
  \end{equation*}
  has codimension $1$ in $U$.
\end{Proposition}

For $U\in \cG_k\cP_n$ and $b\in \Cset$, we inductively def\/ine a f\/lag
of subspaces
\begin{equation*}
  U=U_1 \supset U_2 \supset \cdots \supset U_k,\qquad \dim U_i = k+1-i,
\end{equation*}
by setting $U_1=U$ and
\begin{equation}
  \label{eq:Uidef}
  U_{i+1} = U_i',\qquad i=1,\ldots, k-1.
\end{equation}
Set
\begin{gather*}
  \nu_i  = \ord_b U_i,\qquad
  \lambda_i  = n-\nu_i+i-k,\qquad i=1,\ldots, k.
\end{gather*}
Observe that $\bnu\in \P_{k,\ell}$, that $\blambda\in \B_{k,\ell}$ and
that $\bnu=\blambda^+$,  c.f.~\eqref{eq:lambda+def}.  For reasons
explained below, we will call $\bnu$ the $b$-\emph{pivot sequence} and
call $\blambda$ the $b$-\emph{shape} of $U$.

For polynomials $p_1,\ldots, p_k$, let
\begin{equation}
  \label{eq:Mp1pkdef}
  \M_n(p_1,\ldots, p_k) =
  \begin{pmatrix}
    p_1 & p_1' & \cdots & p^{(n)}_1 \\
    p_2 & p_2' & \cdots & p^{(n)}_2 \\
    \vdots & \vdots & \ddots & \vdots \\
    p_k & p_k' & \cdots & p^{(n)}_k
  \end{pmatrix}
\end{equation}
denote the $k\times (n+1)$ matrix of indicated derivatives.  For $U\in
\cG_k\cP_n$ and $b\in \Cset$, let $\{U_i\}$ be the $\ord_b$ f\/iltration
def\/ined by~\eqref{eq:Uidef}.

\begin{Proposition}
  \label{prop:redbasis}
  There exists a unique basis $p_1,\ldots, p_k$ of $U$ such that
  $p_i\in U_i$ and such that
  \begin{gather*}
    p_i^{(\nu_i)}(b)  = 1,\qquad
    p_i^{(\nu_j)}(b)  = 0,\qquad j> i.
  \end{gather*}
  This basis is characterized by the property that $\M_n(p_1,\ldots,
  p_k)(b)$ is in RREF.
\end{Proposition}

\begin{Definition}
  Henceforth, we will refer to the above $p_1,\ldots, p_k$ as the
  $b$-\emph{reduced basis} of $U$ and set
  \begin{equation*}
    \M(U,b)=\M_n(p_1,\ldots,p_k)(b)\in \R_{k,\ell}.
  \end{equation*}
\end{Definition}

The sequence $\{\nu_i\}$ enumerates the pivot columns of $\M(U,b)$.
Also observe that if we delete the pivot columns from $\M(U,b)$, we
obtain a $\blambda$-bounded matrix.

\begin{Proposition}
  \label{prop:redbasis2}
  Fix $b\in \Cset$. Let $\blambda$ be the $b$-shape of $U$, and
  $p_1,\ldots, p_k$ the $b$-reduced basis. Then,
  \begin{equation}
    \label{eq:piaijdef}
    p_i(z) = \frac{(z-b)^{\nu_i}}{\nu_i!}+ \sum_{j=1}^{\ell} a_{ij}
    \frac{(z-b)^{\mu_j}}{\mu_j!},\qquad i=1,\ldots, k.
  \end{equation}
  where
  $\bnu=\blambda^+$,  $\bmu=\blambda^{\!*\,-}$, and  where the matrix
  \begin{equation}
    \label{eq:aijdef}
    a_{ij} =  p_i^{(\mu_j)}(b),\qquad
    i=1,\ldots, k,\quad j=1,\ldots,\ell,
  \end{equation}
  is  $\blambda$-bounded.
\end{Proposition}

\begin{proof}
  By Theorem \ref{thm:conjpart}, $\bmu$ enumerates the non-pivot
  columns of $M(U,b)$.
  Relation \eqref{eq:piaijdef} follows by \eqref{eq:aijdef} and by
  Proposition \ref{prop:redbasis}.    Observe that, by \eqref{eq:lambda+def}
  \begin{gather*}
    \lambda_i  = n-\nu_i -(k-i) = \#\text{(columns to the right of $\nu_i$)} - \# \text{(pivots to
      the right of $\nu_i$)}\\
   \phantom{\lambda_i}{} = \#\text{(non-pivots to the right of $\nu_i$)}
  \end{gather*}
  More formally,
  \[ \lambda_i = \# \{j: \mu_j> \nu_i \}.\]
  Hence $\mu_j\leq \nu_i$ if and only if $j\leq \ell -\lambda_i$.
  However, $\mu_j$ cannot possibly equal $\nu_i$. Hence, if $j\leq
  \ell -\lambda_i$, then $\mu_j < \nu_i=\ord_b p_i$, and therefore,
  $a_{ij} = p_i^{(\mu_j)}(b) = 0$.
\end{proof}

\begin{theorem}
  \label{thm:dual}
  Let $U\in \cG_k\cP_n$ and $U^*\in \cG_{\ell}\cP_n$ be dual
  polynomial subspaces.  Fix $b\in \Cset$ and let $\blambda$ be the
  $b$-shape of $U$.  Then, the conjugate partition $\blambda^*$ is the
  $b$-shape of the apolar dual $U^*$.  Furthermore, with
  $\bnu$, $\bmu$, $\{ a_{ij}\}$  as above, and with
  $\bnu^*=\blambda^{*+}$, $\bmu^*=\blambda^{-}$ and
  \begin{equation}
    \label{eq:aa*rel}
    a^*_{k+1-i,\ell+1-j} = (-1)^{\nu_i+\mu_j+1} a_{ij},
  \end{equation}
  the reduced basis of $U^*$ is given by
  \begin{equation}
    \label{eq:p*jdef}
    p^*_j(z) = \frac{(z-b)^{\nu^*_j}}{\nu^*_j!} + \sum_{i=1}^{k}
    a^*_{ij} \frac{(z-b)^{\mu^*_i}}{\mu^*_i!},\qquad j=1,\ldots,\ell.
  \end{equation}
\end{theorem}

\begin{proof}
  Since the bilinear form $\gamma$ described in
  \eqref{eq:gammabracket} is $\SL_2\Cset$ invariant, we also have
  \begin{equation}
    \label{eq:gammabracketshifted}
    n!\, \gamma\left(\frac{(z-b)^j}{j!},\frac{(z-b)^k}{k!}\right)= \begin{cases}
      \displaystyle (-1)^j,  & \text{if} \ \ j+k=n,\\
      0,  & \text{otherwise.}
    \end{cases}
  \end{equation}
  Recall that
  \[ \nu_i +\mu^*_{k+1-i}=n,\qquad \mu_j+\nu^*_{\ell+1-j}=n.\] Also, by
  Theorem~\ref{thm:conjpart}, $\{\nu^*_j\}$ and $\{ \mu^*_i\}$ are
  complementary enumerations of the set $\{ 0,1,\ldots, n\}$. Hence
  for f\/ixed $i$, $j$ we have, by the def\/initions \eqref{eq:piaijdef},
  \eqref{eq:p*jdef} and by~\eqref{eq:gammabracketshifted},
  \begin{gather*}
    \gamma(p_i,p^*_{\ell+1-j})  = (-1)^{\nu_i} a^*_{k+1-i,\ell+1-j} +
    (-1)^{\mu_{j}} a_{i,j} =0.\tag*{\qed}
  \end{gather*}
  \renewcommand{\qed}{}
\end{proof}

The above discussion, based on a f\/inite parameter $b\in \Cset$, has an
analogous development in terms of $\infty$.  Fix $k$, $\ell$, $n$ as above.
For a given polynomial $p\in \cP_n$, we call $\infty$ a root of
multiplicity $\nu$ if $\deg p = n-\nu$.  In this sense,
\begin{equation}
  \label{eq:ordinfty}
  \ord_\infty p = n-\deg p.
\end{equation}
Let us also def\/ine
\begin{equation*}
  \deg U = \max \{ \deg p: p \in U \}.
\end{equation*}
Equivalently, if $U$ is non-trivial, then $m=\deg U$ is the smallest
integer such that $U\subset \cP_m$.  We come to the following analogue
of Proposition~\ref{prop:ordU}.
\begin{Proposition}
  Let $U$ be a finite-dimensional polynomial subspace. Then, the
  subspace
  \[ U'=\{ p\in U : \deg p < \deg U \}\]
  has codimension~$1$ in~$U$.
\end{Proposition}

For $U\in \cG_k\cP_n$, we def\/ine a f\/lag of subspaces $\hU_k\subset
\cdots \subset \hU_1$ by setting $\hU_1=U$ and
\begin{equation*}
  \hU_{i+1} = \hU'_i,\qquad i=1,\ldots, k-1.
\end{equation*}
We def\/ine the the $\infty$-shape and the $\infty$-pivot sequence as
follows:
\begin{gather}
   \label{eq:lambdaidegdef}
   \hlambda_i   =\deg \hU_i-k+i,\qquad i=1,\ldots, k,\\
  \hnu_i     =   \deg \hU_{k+1-i}.\nonumber
\end{gather}
Observe that $\hat{\bnu} = \hat{\blambda}{}^{\!-}$.
\begin{Proposition}
  There exists a unique basis $p_1,\ldots, p_k\in U$ such that each
  $p_i\in \hU_{k+1-i}$ and such that
  \begin{gather*}
     p_i^{(\hnu_i)}(0) = 1,\qquad
     p_i^{(\hnu_j)}(0) = 0,\qquad j<i.
  \end{gather*}
\end{Proposition}

We will refer to such a basis as the $\infty$-\emph{reduced
  basis} of $U$.
We can now state the
following analogue of Theorem~\ref{thm:dual}.
\begin{theorem}  \label{thm:dual1}
  Let $U\in \cG_k\cP_n$ and $U^*\in \cG_{\ell}\cP_n$ be dual
  polynomial subspaces.  Then $\hat{\blambda}$ and $\hat{\blambda}^*$,
  the respective $\infty$-shapes, are conjugate partitions.
  The respective $\infty$-reduced
  bases are given by
  \begin{gather}
    \label{eq:pilead}
    p_i(z)  = \frac{z^{\hnu_i}}{\hnu_i!}+\sum_{j=1}^{\ell}
    a_{ij} \frac{z^{\hmu_j}}{\hmu_j!}, \qquad i=1,\ldots, k,\quad
    \hat{\bnu} = \hat{\blambda}{}^{\!-},\quad \hat{\bmu} =
    \hat{\blambda}{}^{\!*+},\\
    p^*_j(z)  = \frac{z^{\hnu^*_j}}{\hnu^*_j!} + \sum_{i=1}^{k}
    a^*_{ij} \frac{z^{\hmu^*_i}}{\hmu^*_i!},\qquad j=1,\ldots,
    \ell,\quad \hat{\bnu}{}^* = \hat{\blambda}{}^{\!*-},\quad
    \hat{\bmu}^* = \hat{\blambda}{}^{\!+},\nonumber
  \end{gather}
  where $a_{ij}$ and $a^*_{ij}$ are related by \eqref{eq:aa*rel}.
\end{theorem}

Finally, we combine order and degree reductions to obtain a useful
characterization of subspaces generated by monomials.  Let $U\in
\cG_k\cP_n$ be a polynomial subspace.  Let $\nu_1<\cdots < \nu_k$ be
the $0$-pivot sequence and let $\hnu_1<\cdots < \hnu_k$ be the
$\infty$-pivot sequence.  These correspond to the following bases of~$U$:
\begin{gather*}
  p_i(z)  = z^{\nu_i}/\nu_i! + \text{higher order terms}, \\
  \hp_i(z)  = z^{\hnu_i}/\hnu_i! + \text{lower degree terms}.
\end{gather*}

\begin{Proposition}
  \label{prop:monom1}
  For $j=1,\ldots, k$, we have $\nu_j \leq \hnu_{j}$.  If
  $\nu_j=\hnu_j$ for a particular $j$, then the monomial $z^j\in U$ .
\end{Proposition}

\begin{proof}
  For each $j=1,\ldots,k$ choose a $p_j\in U$ such that $\ord_0 p_j =
  \nu_j$ and such that $\pi_j:=\deg p_j$ is as small as possible.
  Observe that for $i< j$ we must have $\pi_i \neq \pi_j$; otherwise a
  linear combination of $p_i$ and $p_j$ would have the same order as
  $p_i$ but a smaller degree.  It follows that $\pi_1,\ldots, \pi_k$
  enumerates the degree pivot set $\{ \hnu_j\}$, although not
  necessarily in any particular order.  Observe that for every
  polynomial $p$, we have $\ord_0 p \leq \deg p$ with equality if and
  only if~$p$ is a~monomial. Hence, $\nu_j \leq \pi_j$ with equality
  if and only if $p_j$ is a monomial.  Hence, for $i\leq j$ we have
  \[ \nu_i \leq \nu_j \leq \pi_j.\] Since $\{ \pi_{i},\ldots, \pi_k\}$ is
  a set with $k+1-i$ distinct elements, by the pigeonhole principle,
  \[ \nu_i \leq \min \{ \pi_i,\ldots, \pi_k \} \leq \tnu_{i},\qquad
  i=1,\ldots, k.\] As for the f\/inal assertion, observe that if $i<j$,
  then $\nu_i < \pi_j$.  Hence, if $\nu_i = \hnu_i$, then necessarily
  $\pi_i = \nu_i$; this means that $p_i$ is a monomial.
\end{proof}
\begin{Corollary}
  \label{cor:monom2}
  A polynomial subspace $U\in \cG_k\cP_n$ is spanned by monomials
  if and only if the order pivots match the  degree pivots:  $\nu_i =
  \hnu_i$ for all $i=1,\ldots, k$.
\end{Corollary}

\section{The Wronskian covariant}
\label{sect:inv}

In this section we introduce the Wronskian of a polynomial subspace, a
key covariant that will allow us to solve the equivalence problem.
Throughout, integers $k$, $\ell$, $n$ are as in the preceding section.  For
an ordered set of polynomials $p_1,\ldots, p_k\in \cP_n$, we def\/ine
their Wronskian to be the polynomial
\begin{equation*}
  \W(p_1,\ldots, p_k) = \det \M_{k-1}(p_1,\ldots, p_k);
\end{equation*}
see \eqref{eq:Mp1pkdef} for the def\/inition of $\M$.
For an  integer sequence $\mu_1,\ldots, \mu_k$ def\/ine
\begin{equation*}
  \V(\bmu)=\V(\mu_1,\ldots, \mu_k)  = \prod_{1\leq i<j\leq k} (\mu_j-\mu_i).
\end{equation*}
Observe that if $\bmu$ is strictly increasing, then $V(\bmu)>0$.

\begin{Proposition}
  \label{prop:monwronsk}
  For a bounded partition $\blambda\in \B_{k,\ell}$ and $\bnu =
  \blambda^-$ we have
  \begin{equation*}
    \W(z^{\nu_1},\ldots, z^{\nu_k}) = \V(\bnu)
    z^{N(\blambda)}.
  \end{equation*}
\end{Proposition}

\begin{proof}
  By def\/inition,
  \begin{equation}
    \label{eq:Wexpand}
    \W(z^{\nu_1},\ldots, z^{\nu_k})=\sum_\pi \prod_{i=1}^k
    (\nu_{\pi_i})_{i-1} z^{\nu_{\pi_i}-i+1},
  \end{equation}
  where
  \begin{equation*}
    (x)_{j} = x(x-1)\cdots (x-j+1)
  \end{equation*}
  is the usual falling factorial, and where $\pi$ ranges over all
  permutations of $\{1,\ldots, k\}$.
  By \eqref{eq:lambda-def},
  \begin{equation*}
    \sum_{i=1}^k (\nu_{\pi_i}-i+1) = \sum_{i=1}^k (\nu_i-i+1) =
    \sum_{i=1}^k \lambda_{k+1-i} = N(\blambda).
  \end{equation*}
  Hence, by inspection,
  \[\W(z^{\nu_1},\ldots, z^{\nu_k}) = P(\nu_1,\ldots, \nu_k) z^N,\]
  where $P$ is a polynomial of degree $k$ in each $\nu_i$.  Since the
  Wronskian is an alternating mapping, the function $P$ must be an
  alternating polynomial.  Hence, $\nu_i-\nu_j$ is a factor of $P$ for all
  $i\neq j$, and  hence, $P$ and $\V(\nu_1,\ldots, \nu_k)$ agree up to a
  constant factor.  By comparing the highest order terms of $\nu_k$ in
  both $\V(\nu_1,\ldots, \nu_k)$ and \eqref{eq:Wexpand} we see that this
  factor is $1$.
\end{proof}
\begin{Proposition}
  \label{prop:degW}
  Let $p_1,\ldots, p_k\in \cP_n$ be linearly independent polynomials,
  and let $\hat{\blambda}$ be the $\infty$-shape of the subspace
  $U=\lspan\{ p_1,\ldots, p_k\}$ as defined in
  \eqref{eq:lambdaidegdef}.  Then,
  \begin{equation}
    \label{eq:degW}
    \deg \W(p_1,\ldots, p_k) =\sum_i \hlambda_i\leq k\ell.
  \end{equation}
  Equality holds if and only if $p_1,\ldots, p_k$ are linearly
  independent modulo  $\cP_{\ell-1}$.
\end{Proposition}
\begin{proof} It suf\/f\/ices to establish \eqref{eq:degW} in the case
  where $p_1,\ldots, p_k$ is the $\infty$-reduced basis.  Let $\hnu_i
  = \deg p_{k+1-i}$ be the pivot sequence.  By the multi-linearity of
  the Wronskian and by Proposition~\ref{prop:monwronsk},
  \[ \W(p_k,\ldots, p_1) = \frac{\V(\hnu_1,\ldots, \hnu_k)}{\hnu_1!\cdots
    \hnu_k!}  z^{N(\hat{\blambda})}+\text{ lower degree terms.} \]
  This establishes the equality in~\eqref{eq:degW}.

  We now prove the second assertion.  Since $\hat{\blambda}$ is a
  $k\times \ell$ bounded partition, we have $N(\hblambda)\leq k\ell$
  with equality if and only if $\hblambda$ is rectangular.  By
  \eqref{eq:pilead} and \eqref{eq:lambda-def}, this is true if and
  only if $\hnu_1=\hlambda_k=\ell$.  Since $\hnu_1<\cdots <\hnu_k \leq
  n$ and since $\ell+k-1 =n$, the partition $\hblambda$ is rectangular
  if and only if $p_1,\ldots, p_k$ are linearly independent modulo
  $\cP_{\ell-1}$.
\end{proof}

Let $U\in\cG_k\cP_n$ be a polynomial subspace and let $p_1,\ldots,
p_k\in U$ be a basis.  In light of Proposition~\ref{prop:degW},
we may regard the Wronskian $\W(p_1,\ldots, p_k)$ as an element of
$\cP_{k\ell}$.  Since $\W$ is a multi-linear, alternating mapping,
\begin{equation*}
  p_1\wedge \cdots \wedge p_k \mapsto \W(p_1,\ldots, p_k),\qquad
  p_1,\ldots,p_k\in \cP_n
\end{equation*}
def\/ines a linear transformation $\W:\Lambda^k\cP_n\to \cP_{k\ell}$;
here $\Lambda^k$ denotes the $k\supth$ exterior product of a~vector
space.

Let us give a basis description of this linear transformation.  For a
sequence of integers $0\leq \nu_1,\ldots ,\nu_k\leq n$ def\/ine
\begin{equation*}
  z^\bnu = z^{\nu_1} \wedge \cdots \wedge z^{\nu_k}\in \Lambda^k\cP_n.
\end{equation*}
In light of the bijective correspondence $\B_{k,\ell}\to
\P_{k,\ell}$ def\/ined in Section \ref{sub:partitions}, the set $\{
z^{\blambda^-} : \blambda \in \B_{k,\ell} \}$ is a basis of
$\Lambda^k\cP_n$.  We can therefore describe $\W$ by stipulating
\begin{equation*}
  \W(z^{\blambda^-}) = \V(\blambda^-)  z^{N(\blambda)},\qquad
  \blambda\in \B_{k,\ell}.
\end{equation*}

\begin{theorem}
  \label{thm:Wequivariant}
  The linear transformation $\W:\Lambda^k \cP_n \to \cP_{k\ell}$ is
  $\SL_2\Cset$-equivariant.
\end{theorem}

\begin{proof}
  Let
  \begin{equation}
    \label{eq:fracxform}
    z = \frac{a\hz+b}{c\hz + d},\qquad ad-bc=1,
  \end{equation}
  be a fractional linear transformation, let $p_1,\ldots, p_k\in
  \cP_n$, and let
  \[ \hp_i(\hz) = (c\hz+d)^n p_i(z),\qquad i=1,\ldots, k\] be the
  transformed polynomials.  The claim is that
  \begin{equation}
    \label{eq:Wequiv}
    \W(\hp_1,\ldots, \hp_k)(\hz) = (c\hz+d)^{k\ell} \W(p_1,\ldots,
    p_k)(z).
  \end{equation}

  It suf\/f\/ices to consider 3 kinds of transformation: translations,
  homotheties, and inversions.  Consider a translation transformation
  \[ \hp(\hz) = p(z),\qquad p\in \cP_n,\qquad z = \hz+b,\qquad b\in
  \Cset.\] Here $a=1$, $c=0$, $d=1$.  Observe that $\hp^{(j)}(\hz) =
  p^{(j)}(z)$.  Equation~\eqref{eq:Wequiv} follows immediately.

  Next, consider a homothety transformation
  \begin{equation*}
     \hp(\hz) = a^{-n} p(z),\qquad z= a^2 \hz,\qquad a\neq 0.
  \end{equation*}
  This corresponds to $b=c=0$ and $d=1/a$ in
  \eqref{eq:fracxform}. Observe that \[\hp^{(j)}(\hz) = a^{2j-n}
  p^{(j)}(z).\] Hence,
  \[ \W(\hp_1,\ldots, \hp_k)(\hz) = a^{k(k-1)-nk} \W(p_1,\ldots,
  p_k)(z) = a^{-k\ell} \W(p_1,\ldots, p_k)(z),\] in full  accordance
  with~\eqref{eq:Wequiv}.

  Finally, consider an inversion
  \[ \hp(\hz) =\hz^n p(z),\qquad z=-1/\hz.\] Let $R_i(z)$, $
  i=0,1,\ldots, k-1$ denote the $i\supth$ column vector
  $(p^{(i)}_1(z),\ldots, p^{(i)}_k(z))^t$ of $\M_{k-1}(p_1$, $\ldots,
  p_k)$, and let $\hR_i(\hz)$ be the corresponding column vector formed
  from the derivatives of $\hp_i(\hz)$.  Observe that
  \[ \hp'(\hz) = \hz^{n-2} p'(z) + n \hz^{n-1} p(z),\]
  whence
  \[ \hR_1(\hz) = \hz^{n-2} R_1(z) + n \hz^{n-1} R_0(z).\]
  More generally,
  \[ \hR_i(\hz) = \hz^{n-2i} R_i(z) + \cdots ,\qquad  i =0,\ldots, k-1,\]
  where the remainder is a linear combination of $R_0(z),\ldots,
  R_{i-1}(z)$.  Since
  \[ \sum_{i=0}^{k-1} n-2i = nk - k(k-1) = k\ell,\] by the
  multi-linearity and skew-symmetry of the determinant it follows that
  \begin{gather*}
    \W(\hp_1,\ldots, \hp_k)(\hz)  =
    \det(\hR_0(\hz),\ldots,\hR_{k-1}(\hz)) = \hz^{k\ell} \W(p_1,\ldots, p_k)(z),
  \end{gather*}
  as was to be shown.
\end{proof}

\begin{Definition}
  By a slight abuse of notation we also use $\W$ to denote the
  corresponding Wronski mapping $\W:\cG_k\cP_n \to \cG_1 \cP_{k\ell}$,
  def\/ined by
  \[ \W(U) = \lspan \W(p_1,\ldots, p_k),\qquad U\in \cG_k\cP_n,\] where
  $p_1,\ldots, p_k\in U$ is any basis.  We refer to $\W(U)$ as the
  \emph{Wronskian covariant} of the polynomial subspace $U$.
\end{Definition}

Being an element of $\cG_1\cP_{k\ell}$, we may regard $\W(U)$ as a
polynomial of degree $\leq k\ell$ def\/ined up to a non-zero scalar
multiple.  However, the $k\ell$ roots of $\W(U)$ are def\/ined
unambiguously.  These roots transform covariantly with respect to
projective transformations.

The relationship of the linear and non-linear Wronski maps is
conveniently summarized by the following diagram; the top arrow is the
Pl\"ucker embedding, the right arrow is the projectivization of the
linear Wronski map and the left diagonal arrow is the non-linear
Wronski map
\[\xymatrix{
  \cG_k\cP_n \ar @{^{(}->}[r]
  \ar[dr]
  & \cG_1(\Lambda^k\cP_n) \ar[d]\\
  & \cG_1\cP_{k\ell}}
  \]

We can now formulate the relationship between roots of the Wronskian
and the  possible shapes of a subspace.  Throughout, recall that the
order of a polynomial at $\infty$ is def\/ined by~\eqref{eq:ordinfty}.

\begin{theorem}
  \label{thm:Word}
  Let $U\in \cG_k\cP_n$ be a subspace and let $b\in \Cset \cup \{
  \infty \}$ be given.  We have
  \begin{equation*}
    \ord_{b} \W(U) = N(\blambda') =  k\ell - N(\blambda),
  \end{equation*}
  where $\blambda$ is the $b$-shape of $U$ and $\blambda'$ is the
  complementary partition.
\end{theorem}
\begin{proof}
  This follows from Proposition \ref{prop:degW}, which treats the
  $b=\infty$ case.  By Proposition \ref{prop:degW},
  \[ \ord_\infty \W(U) = k\ell-N(\blambda),\] where $\blambda$ is the
  $\infty$-shape of $U$.  The present assertion now follows by the
  equivariance of $\W$.
\end{proof}

The description of simple roots deserves explicit mention.
\begin{Corollary}
  Let $U\in\cG_k\cP_n$ be a subspace. Fix $b\in \Cset$ and let
  $\blambda$ be the $b$-shape of $U$ and $\bnu=\blambda^+$ the pivot
  sequence. Then, $b$ is a simple root of $\W(U)$ if and only if
  $\nu_j = j-1$, $j=1,\ldots, k-1$ and $\nu_k = k$.  Equivalently, $b$
  is a simple root if and only if $\M(U,b)$ has the form
  \begin{equation*}
    \M(U,b) =
    \begin{pmatrix}
      I_1 & B_{11} & 0 & B_{12} \\
      0 & 0 &  1 & B_{22}
    \end{pmatrix},
  \end{equation*}
  with $B_{11}$, $B_{12}$, $B_{22}$ having dimensions $(k-1)\times 1$,
  $(k-1)\times (\ell-1)$ and $1\times (\ell-1)$, respectively; c.f.~\eqref{eq:RREF}.
\end{Corollary}

We also have the following key result asserting the invariance of the
Wronskian with respect to apolar duality.

\begin{theorem}
  \label{thm:Wdual}
  A polynomial subspace $U\in\cG_k\cP_n$ and its apolar dual $U^*\in
  \cG_{\ell}\cP_n$ have the same Wronskian covariant; $\W(U) =
  \W(U^*)$.
\end{theorem}

\begin{proof}
  Theorems \ref{thm:dual} and \ref{thm:dual1} establish shape duality
  of $U$ and $U^*$.  By Proposition \ref{prop:conjpartchar}, we have
  $N(\blambda)=N(\blambda^*)$. Therefore, by Theorem \ref{thm:Word}
  the Wronskians $\W(U)$ and $\W(U^*)$ have the same roots with the
  same multiplicities.
\end{proof}

\subsection{Additional properties of the Wronski map}
Call a multi-vector in $\Lambda^k\cP_n$ \emph{decomposable} if it can
be given as $p_{1} \wedge \cdots \wedge p_k$ for some linearly
independent $p_1,\ldots, p_k \in \cP_n$.  A basis of a $k$-dimensional
subspace $U\in \cG_k \cP_n$ def\/ines a decomposable multi-vector,
unique up to scalar multiple.  This gives us the Pl\"ucker embedding
$\cG_k \cP_n\hookrightarrow\cG_1 \Lambda^k \cP_n$.  Indeed, the
Grassmannian $\cG_k \cP_n$ can be regarded as a projective variety
generated by a collection of quadratic polynomials, the so-called
Pl\"ucker relations \cite{eastwood,griffiths}.

The domain and codomain of the Wronski mapping $\W :\cG_k\cP_n \to
\cG_1 \cP_{k\ell}$ are both $k\ell$-dimensional varieties.  This
suggests that $\W$ has nice properties from the algebro-geometric
point of view. The following result is needed in the sequel.
\begin{Proposition}
  \label{prop:wronsksurj}
  The Wronski map $\W:\cG_k \cP_n \to \cG_1 \cP_{k\ell}$ is
  surjective with finite pre-images of points.
\end{Proposition}

\begin{proof}
  The kernel of the linear transformation $\W:\Lambda^k\cP_n \to
  \cP_{k\ell}$ does not contain any pri\-mi\-tive multi-vectors; see
  Proposition \ref{prop:monwronsk}. Hence, the mapping $\W:\cG_k\cP_n
  \to \cG_1\cP_{k\ell}$ is a central projection, and hence is f\/inite
  and surjective~\cite[Theorems 4 and 7, Section 5.3]{shafarevich}.
\end{proof}

Consideration of the degree of the Wronski map leads to very
interesting topics such as the combinatorics of Young tableaux and the
Schubert calculus.  Indeed, the degree is given by $d(k,l)$ where
\begin{equation}
  \label{eq:dkldef}
  d(k,\ell) = \frac{(k-1)!! (\ell-1)!!}{n!!} (k\ell)!,\qquad\text{where} \quad j!! = 1!  2! \cdots j!,
\end{equation}
counts the standard Young tableaux of shape $k\times \ell$. See
references \cite{fulton,griffiths} for more details.

Finally, let us mention the following, alternate characterization of
the Wronski map \cite{eremenko}.
\begin{Proposition}
\label{prop:altwronsk}
For $p_1,\ldots, p_k \in \cP_n$, we have
  \begin{gather}
    W(p_1,\ldots, p_k)(z) (w^0\wedge w^1\wedge \cdots \wedge w^n) \nonumber
    \\
    {} \qquad = (k-1)!! \, p_1(w)
    \wedge \cdots  \wedge p_k(w) \wedge (w-z)^k\wedge \cdots
    \wedge (w-z)^n.\label{eq:altwronsk}
  \end{gather}
\end{Proposition}

\begin{proof}
  Note that
  \[ (w-z)^j = \sum_{i=0}^j \binom{j}{i} (-1)^i w^{j-i} z^{i} \] is
  understood to refer to an element of $\Lambda^1(\cP_n(w))\otimes \cP_n(z)$,
  while relation \eqref{eq:altwronsk} should be understood to hold in
  $\Lambda^k\left(\cP_n(w)\right) \otimes \cP_n(z)$. In other words,
  the wedge product is taken relative to the $w$ variable.  For
  example:
  \[ (aw+bz) \wedge (w-z)^2 = a (w^1 \wedge w^2) + bz (w^0 \wedge w^2) - (a+2b)
  z^2 (w^0\wedge w^1).\]
  Relation \eqref{eq:altwronsk} follows from the following observations:
  \begin{gather*}
     p_i(w) = \sum_{j=0}^{k-1} p_i^{(j)}(z) \frac{(w-z)^j}{j!},\qquad
    i=1,\ldots, k ,\\
     w^0\wedge w^1 \wedge \cdots \wedge w^n = (w-z)^0\wedge (w-z)^1 \wedge
    \cdots \wedge (w-z)^n .\tag*{\qed}
  \end{gather*}
  \renewcommand{\qed}{}
\end{proof}

Formula \eqref{eq:altwronsk} has the following, geometric
interpretation.  Let us regard $\cG_1\cP_n$ as $n$-di\-men\-sio\-nal
projective space and via projectivization identify $k$-dimensional
subspaces $U\in \cG_k \cP_n$ with $(k-1)$ dimensional f\/lats in
$\cG_1\cP_n$.  For a primitive multivector $p_1\wedge\cdots\wedge
p_k$, let $[p_1\wedge \cdots \wedge p_k]$ denote the corresponding
element of $\cG_k\cP_n$ under the Pl\"ucker embedding. The mapping
\[ z\mapsto [(w-z)^n] \] describes the rational normal curve in
projective space $\cG_1\cP_n(w)$.  Setting \[ N(z) = (w-z)^k\wedge
\cdots \wedge (w-z)^n\in \Lambda^\ell \cP_n(w),\] the $1$-parameter
family of corresponding $(\ell-1)$ dimensional osculating f\/lats is
described by $[N(z)]\in \cG_\ell\cP_n(w)$.  Therefore, relation
\eqref{eq:altwronsk} implies that the roots of the Wronskian
$W(U)(z)$, where $U\in \cG_k\cP_n$, correspond to points on the
rational normal curve where the $(\ell-1)$ dimensional f\/lat $[N(z)]$
touches the $(k-1)$ dimensional f\/lat $U$ in the ambient $(k+\ell-1)$
dimensional projective space $\cG_1\cP_n$.  The order of the root
corresponds to the dimension of the intersection plus~1.  Finally,
note that since the rational normal curve is $\SL_2\Cset$-invariant,
the above observations constitute an alternate proof of the
equivariance of the Wronski map~$\W$.

\section{The equivalence problem}
\label{sect:eqprob}
\subsection{The generic reduction}
\begin{Definition}
  Let $U,\hU\in \cG_k\cP_n$ be polynomial subspaces.  Recall that
  $U\sim \hU$ means that there exists a transformation $g\in
  \PSL_2\Cset$ such that $\hU = g\cdot U$.  We def\/ine $\Stab(U)\subset
  \PSL_2\Cset$ to be the stabilizer of $U$; that is, the group of
  transformations $g\in \PSL_2\Cset$ such that $U=g\cdot U$.
\end{Definition}

In light of Proposition \ref{prop:wronsksurj}, every polynomial $Q\in
\cP_{k\ell}(z)$ corresponds to a f\/inite number of subspaces $U\in
\cG_k\cP_n$ such that $\W(U) =\lspan Q$.  In particular, if $Q$ has
simple roots, then the number of such subspaces, counting
multiplicities, is $d(k,\ell)$ as given by
\eqref{eq:dkldef}\footnote{For some values of the roots of W(z) the
  number of distinct elements may decrease. In fact it is an
  interesting question to understand when the roots are non-generic.
  We thank the referee for this observation.}.  The equivariance of the
Wronski map, Theorem \ref{thm:Wequivariant}, has the following
immediate consequences.
\begin{Proposition}
  For $U\in \cG_k\cP_n$ we have $\Stab(U) \subset \Stab(W(U))$.
\end{Proposition}

The above observation leads to an immediate proof of Proposition
\ref{prop:genequiv}.  Since a generic polynomial has a trivial
symmetry group, the algorithm of~\cite{BO00} reduces the
equivalence problem for generic polynomial subspaces to the
corresponding equivalence problem for the corresponding Wronskians.
Again, it must be emphasized that if $W(U)$ is equivalent to $W(\hU)$,
where $U,\hU\in \cG_k\cP_n$, we cannot conclude that $U\sim \hU$.
Rather, in the generic case it suf\/f\/ices to calculate the unique $g\in
\PSL_2$ such that $W(\hU) = g\cdot W(U)$ and then to test whether $\hU
= g\cdot U$.

Let us recall the relevant details of the algorithm of
\cite{BO00,olver99}.  For $Q\in \cP_N$\footnote{In the present setting
  we are interested in the case where $N=k\ell$ is the degree of the
  Wronskian covariant.}, the following dif\/ferential operators def\/ine a
number of key polynomial covariants \cite[Chapter 5]{olver99}:
\begin{gather}
  \label{eq:Hdef}
  H  = \frac{1}{2}(Q,Q)^{(2)}
   = N(N-1)\left[QQ''-\frac{N-1}{N}
    (Q')^2\right],\\
  T  = (Q,H)^{(1)}
   = -N^2(N-1)  \left[ Q^2 Q''' - 3 \frac{(N-2)}{N} Q
    Q' Q'' + 2
    \frac{(N-1)(N-2)}{N^2}   (Q')^3\right],\!\! \\
  V  = Q^3 Q''' - 4 \frac{(N-3)}{N}  Q^2 Q' Q'' + 6
  \frac{(N-2)(N-3)}{N^2}\, Q (Q')^2 Q'' \\ \nonumber
\phantom{V  =}{}  + \frac{(N-1)(N-2)(N-3)}{N^3} (Q')^4,\\
  \label{eq:Udef}
  U = (Q,T)^{(1)}  = N^3(N-1) V- 3\frac{(N-2)}{(N-1)} H^2.
\end{gather}
Above, for polynomials $Q\in \cP_n$ and $R\in \cP_m$ the expression
\begin{gather*}
  (Q,R)^{(r)}(z)   = r! \sum_{k=0}^r (-1)^k \binom{n-r+k}{k}
  \binom{m-k}{r-k} Q^{(r-k)}(z) R^{(k)}(z)
\end{gather*}
is a joint covariant of $Q$ and $R$ called the $r$th
transvectant~\cite[Chapter~5]{olver99}.

We also form the following, fundamental, zero-weight rational covariants:
\begin{gather*}
  J^2  = \frac{T^2}{H^3},\qquad
  K  = \frac{U}{H^2}
\end{gather*}
  Let $Q, \hat{Q}\in \cP_n$ be two polynomials and let $H=H[Q]$,
  $J=J[Q]$, $K=K[Q]$ and $\hat{H} = H[\hat{Q}]$, $\hat{J} = J[\hat{Q}]$,
  $\hat{K} = K[\hat{Q}]$ be the indicated covariants.

  \begin{theorem}[Theorems~2.5 and~5.8 of \cite{BO00}]
    The polynomials $Q$, $\hat{Q}$ are $\PSL_2\Cset$-equivalent if and
    only if both $H,\hat{H}=0$ or both $H,\hat{H}\neq 0$ and the
    following rational relations are compatible
    \begin{gather}
      \label{eq:JKcompat}  J^2(z) = \hat{J}^2(\hz)\qquad\text{and}\qquad
      K(z) = \hat{K}(\hz).
    \end{gather}
    Furthermore, if \eqref{eq:JKcompat} are consistent, then the
    resulting relation between $z$ and $\hz$ is a fractional linear
    transformation.
\end{theorem}

Applying the criterion to $Q=W(U)$ and $\hat{Q} = W(\hU)$ with
$N=k\ell = k (n+1-k)$ provides necessary condition for the equivalence
of two polynomial subspaces $U,\hat{U}\in \cG_k\cP_n$.  Indeed if
$\Stab Q$ is trivial, the criterion is also suf\/f\/icient.  If $\Stab Q$ is
a f\/inite group, then the equivalence problem is reduced to the
following f\/inite procedure.

\begin{theorem}
  \label{thm:genequiv}
  Let $U,\hU\in \cG_k\cP_n$ be polynomial subspaces.  Suppose that
  $W(\hU) = g_1\cdot W(U)$ for some $g_1\in \PSL_2\Cset$.  Then $U\sim
  \hU$ if and only if there exists a projective transformation
  $g_2\in\Stab W(U)$ such that $\hU = g_1\cdot g_2\cdot U$.
\end{theorem}

We refer to \cite{BO00} for the classif\/ication of binary forms
(polynomials) with f\/inite, but non-trivial symmetry groups.

\subsection{The case of inf\/inite symmetries}

In light of the above remarks, a full solution of the equivalence problem
requires an understanding the class of subspaces $U\in \cG_k \cP_n$
with an inf\/inite $\Stab U$.  In this regard, let us recall the
following helpful results from \cite{BO00}.  Throughout, $U\in
\cG_k\cP_n$ is a polynomial subspace, $Q=W(U)$ is its Wronskian
covariant, def\/ined up to scalar multiple, and $N=k\ell=k(n+1-k)$.
\begin{Proposition}
  \label{prop:aut1}
  The following are equivalent:
  \begin{itemize}\itemsep=0pt
  \item $Q(z)$ is a constant or $Q(z) = (az+b)^{N}$ is a monomial
    of maximum degree.
  \item $H[Q] = 0$;
  \item $\Stab(Q)$ is a $2$-dimensional subgroup consisting of fractional
    linear transformations that fix exactly one point in $\hC$.
  \end{itemize}
\end{Proposition}

\begin{Proposition}
  \label{prop:aut2}
  The following are equivalent:
  \begin{itemize}\itemsep=0pt
  \item $Q(z)=(az+b)^m$ is a monomial where $m\neq 0, N$;
  \item $J^2[Q]$ is a constant;
  \item $\Stab(Q)$ is a $1$-dimensional subgroup consisting of fractional linear
    transformations that fix exactly two  points in $\hC$.
  \end{itemize}
\end{Proposition}

 We use the above results to establish the following.
Let us say that $U\in \cG_k\cP_n$ is a monomial subspace if $U =
\lspan\{ z^{\nu_1}, \ldots, z^{\nu_k}\}$, where $0\leq \nu_1<\cdots <
\nu_k \leq n$.
\begin{theorem}
  \label{thm:onetworoots}
  The following are equivalent:
  \begin{itemize}\itemsep=0pt
  \item $U$ is projectively equivalent to a monomial subspace;
  \item $Q=\W(U)$ has either one or two distinct roots;
  \item $H[Q]=0$ or  $J^2[Q]$ is constant;
  \item $\Stab(Q)$ is infinite $($of the type indicated in the above two
    propositions$)$.
  \end{itemize}
\end{theorem}

\begin{proof}
  We focus our attention on the most interesting of the above
  implications.  First, suppose that $U$ is projectively equivalent to
  a monomial subspace.  Since all of the above properties are $\SL_2$
  equivariant, no generality is lost if we assume that $U$ is a
  monomial subspace with pivot sequence $\{\nu_i\}$.  By Proposition
  \ref{prop:monwronsk}, the Wronskian
  \[ \W(U)(z)= [z^m],\qquad\text{where} \quad m = \sum_i \nu_i - k(k-1)/2\]
  is also a monomial. Hence, $0$ and $\infty$ are the only roots.
  The other conclusions follow by Propositions~\ref{prop:aut1} and~\ref{prop:aut2}.

  Next, suppose that $\W(U)$ has one or two distinct roots.  We will
  show that $U$ is a monomial subspace. By virtue of
  $\SL_2$-equivariance, no generality is lost in supposing that these
  roots are~$0$ and~$\infty$ or just $\infty$ if there is only 1 root.
  Let $\{ \nu_i \}$ and $\{\hnu_i\}$ be the corresponding~$0$ and
  $\infty$-pivot sequences.  By Proposition \ref{prop:monom1}, $\nu_i
  \leq \hnu_i$.  Let
 \begin{gather*}
   \lambda_i  = n-\nu_i-(k-i),\qquad
   \hlambda_i  = \hnu_{k+1-i}-(k-i)
 \end{gather*}
 be the corresponding shape partitions.  By Theorem \ref{thm:Word},
 \begin{gather*}
    \ord_0 Q   =  k\ell  - \sum_i \lambda_i,\qquad
    \deg Q  =  \sum_i \hlambda_i,\qquad
    \deg Q - \ord_0 Q  = \sum_i (\hnu_i-\nu_i)
 \end{gather*}
 Therefore, $\deg Q = \ord_0 Q$ if and only if $\nu_i = \hnu_i$ for
 every $i$.  By Corollary~\ref{cor:monom2}, the latter is true if and
 only if $U$ is spanned by monomials.
\end{proof}

We emphasize that the above criterion is fully algorithmic.  Given a
basis of a polynomial subspace $U\in \cG_k\cP_n$, the algorithm
requires us to calculate the Wronskian, and then calculate a~f\/inite
number of dif\/ferential invariants shown in \eqref{eq:Hdef}--\eqref{eq:Udef}.  If $H=0$, then either $U = \cP_{k-1}$ or there
exists a $b\in \Cset$ such that
\[ U = \lspan\{(z-b)^\ell,\ldots, (z-b)^{n}\}.\] If $H\neq 0$, but
$T^2$ is a constant multiple of $H^3$, then $U$ is projectively
equivalent to a monomial subspace.

\subsection{Primitive subspaces}

We now consider two basic questions related to subspace equivalence.
For $m<n$ we have the natural embedding  $\cP_m\subset \cP_n$.  It is
natural to ask whether $\ell=n+1-k$ is the ``true codimension'', or
whether it is possible to reduce the degree $n$ and the codimension
$\ell$ by means of a fractional linear transformation.  Dualizing, we
obtain the following question: when does a~subspace $U\subset \cP_n$
contain a non-zero $n$th power $(az+b)^n\in U$.

\begin{Definition}
  We say that $U$ is an \emph{imprimitive} subspace if there exists a
  subspace $\hU\sim U$ such that $\hU \subset \cP_{m}$ for some $m<n$.
  A \emph{primitive} subspace is one for which no such $\hU$ exists.
  We say that $U$ is \emph{strongly primitive} if both $U$ and its
  apolar dual $U^*$ are primitive subspaces.
\end{Definition}

We also have the following characterizations of the above
concepts\footnote{In the interpretation of a polynomial subspace
  $U\in \cG_k\cP_n$ as a rational curve it is natural to regard the
  dual subspace $U^*\in \cG_\ell \cP^n$ as a linear series that
  ``carves out'' the given $U$.  In this interpretation the condition
  of primitivity is equivalent to the condition that the series be
  basepoint-free; see~\cite{khsott03,eishar83} for
  more details.  We thank the referee for this observation.}.

\begin{Proposition}
  \label{prop:npower0}
  A subspace $U\in \cG_k\cP_n$ is imprimitive if and only if $\deg
  U<n$ or if $\ord_b U>0$ for some $b\in \Cset$.
\end{Proposition}

\begin{Definition}
    Let $U\in \cG_k\cP_n$ be a polynomial subspace.   We say that $U$
    contains a top power if $(z-b)^n\in U$ for some $b\in \Cset$ or if
    $1\in U$.
\end{Definition}

\begin{Proposition}
  \label{prop:npower1}
  A polynomial subspace $U\in \cG_k\cP_n$ contains a top power if and
  only if the apolar dual space $U^*$ is imprimitive.
\end{Proposition}

\begin{Proposition}
  A subspace $U\in \cG_k \cP_n$ is strongly primitive if and only if
  $\M(U,b)$ is non-degenerate RREF for all $b\in \hC$. Here,
  \textit{non-degenerate} means having form \eqref{eq:RREF}.
\end{Proposition}

To conclude this subsection, we present the following necessary, but
not suf\/f\/icient, condition for primitivity.
\begin{Proposition}
  Let $U\in \cG_k\cP_n$ be a subspace.  If all roots of $\W(U)$ have
  multiplicities $<k$ then $U$ is primitive.  If all roots have
  multiplicities $<\min(k,\ell)$, where $\ell=n+1-k$ is the
  codimension, then $U$ is strongly primitive.
\end{Proposition}

\begin{proof}
  If $U$ is imprimitive, then there exists a $b\in \hC$ such that
  $\M(U,b)$ is a degenerate RREF matrix of type \eqref{eq:degen1} or~\eqref{eq:degen3}.  The width of the corresponding shape is $<
  \ell$.  Hence, by Theorem~\ref{thm:Word} we have \[\ord_b \W(U) \geq
  k\ell-k(\ell-1) = k.\] If $U^*$ is imprimitive, then by Theorems
  \ref{thm:dual}, \ref{thm:dual1}, for some $b$, the RREF matrix
  $\M(U,b)$ is degenerate of type \eqref{eq:degen2} or
  \eqref{eq:degen3}.  The height of the corresponding shape is $<k$,
  and hence
  \[ \ord_b \W(U) \geq \ell.\]  Therefore, if either $U$ or $U^*$
  are imprimitive, then for all $b\in \hC$ we have
\begin{gather*}
 \ord_b \W(U) \geq \min(k,\ell).\tag*{\qed}
  \end{gather*}\renewcommand{\qed}{}
\end{proof}

\subsection{An example}

For the purposes of illustration let us consider the classif\/ication
and equivalence  of 2-dimensional subspaces of $\cP_3$; i.e.,
$k=2,\ell=2$.   In this case $d(2,2)=2$; this value corresponds to the
classic result that generically, 4 lines  in 3-dimensional space are
simultaneously touched by exactly 2~lines.  Consider two linearly
independent polynomials   $p_1, p_2\in \cP_3$.  Generically, the Wronskian
$W(p_1,p_2)$ will have four distinct roots, call them $b_1$, $b_2$, $b_3$, $b_4$.
These complex numbers are distinguished by the fact that the
polynomial subspace $U=\lspan \{ p_1, p_2\}$ contains exactly  four (up to
scalar multiple)
polynomials with a double root.  These polynomials have the form
\[ q_i(z) = r_i(z) (z-b_i)^2, \qquad r_i \in \cP_2.\]  Now the subspaces
\[L_i=\{ r(z) (z-b_i)^2: r(z)\in \cP_2\}\] represent 4 lines in
3-dimensional projective space $\cG_1\cP_3$.  The subspace $U$ also
represents a line.  The polynomials $q_1,\ldots, q_4$ correspond to the 4
points where the line $[U]$ touches the four lines $[L_1], \ldots ,
[L_4]$.  Since $d(2,2)=2$, in the generic case, there exists exactly
one other subspace $\hU\in \cG_2\cP_3$ with the same Wronskian. The
lines $[U]$, $[\hU]$ are the unique lines that simultaneously touch
$[L_1],\ldots, [L_4]$.

Let us now apply the above covariant analysis to describe the
2-dimensional polynomial subspaces of $\cP_3$ with an inf\/inite
symmetry group.  A generic subspace $U\in \cG_2\cP_3$ can be described
as the span of polynomials
\begin{gather*}
  p_1(z)  = z^3 + a_{11} z + a_{10},\qquad
  p_2(z)  = z^2 + a_{21} z + a_{20}.
\end{gather*}
The corresponding Wronskian, and covariants are given by
\begin{gather*}
  W  = z^4 +2a_{21}z^3 +(3a_{20} -a_{11})z^2 -2a_{10} z+ (a_{11}
  a_{20}-a_{10}a_{21}) \\
 \phantom{W}{} =z^4 + b_3 z^3 + b_2 z^2 + b_1 z + b_0,\\
  H  = 3(8b_2-3b_3^2)z^4 +12(6b_1-b_2 b_3) z^3 + 6(24 b_0 - 2b_2^2 + 3
  b_1 b_3)z^2  \\
 \phantom{H=}{}  +12 ( 6 b_0 b_3-b_1 b_2 )z+ 3(8 b_0 b_2-3b_1^2).
\end{gather*}
A straight-forward calculation shows that $H=0$ if and only if
\[
a_{20} = -a_{11}/3 = -a_{10}/a_{21} = a_{21}^2/4,
\]
or equivalently, if\/f
\[
W = (z+a_{21}/2)^4,
\]
in full  accordance with Proposition  \ref{prop:aut1}.

Let us now describe all spaces projectively equivalent to a span of
monomials.  The case where~$U$ is equivalent to the span of $\{1, z\}$
corresponds to $H=0$ and is described above.  Let us now consider the
other cases: $(i)$ $U\sim \lspan \{ z,z^2\}$; $(ii)$ $U \sim \lspan\{
1,z^3\} $; $(iii)$ $U\sim \lspan\{1,z^2\}$.  The other possible monomial
subspaces are all projectively equivalent to one of these
possibilities.  For cases $(i)$ and $(ii)$ we have $J^2=0$; for case $(iii)$
$J^2=-16/9$.  One can also arrive at the special values $J^2=0$ and
$J^2=-16/9$ by  applying the Gr\"obner basis algorithm to the principal
ideal
generated  by $J^2 H^3 - T^2$ with $J$ as an extra variable.

Hence, for cases $(i)$ and $(ii)$ we must consider the condition $T=0$, where
\begin{gather*}
   T/36=
  (8b_1 - 4b_2b_3 + b_3^3)z^6
  +  2(16b_0 - 4b_2^2 + 2b_1b_3 + b_2b_3^2)z^5 \\
\phantom{T/36=}{}  - 5(4b_1b_2 - 8b_0b_3 - b_1b_3^2)z^4
  - 20(b_1^2 - b_0b_3^2)z^3    - 5(8b_0b_1 + b_1^2b_3 - 4b_0b_2b_3)z^2\\
 \phantom{T/36=}{}
  - 2(16b_0^2 + b_1^2b_2 -  4b_0b_2^2 + 2b_0b_1b_3)z
 -b_1^3 + 4b_0b_1b_2 - 8b_0^2b_3.
\end{gather*}
Straightforward calculations show that the general solution consists
of two possibilities:
\begin{gather*}
   a_{11} = a_{20}-a_{21}^2 ,\qquad a_{10} = -a_{20} a_{21},\qquad
  W=(z^2+a_{21}z + a_{20})^2, \\
   a_{11} = 3(3a_{20}-a_{21}^2),\qquad a_{10} = 3
  a_{20}a_{21}-a_{21}^3,\qquad W=(z^2+a_{21}z-a_{11}/3)^2.
\end{gather*}
Letting $r_1$, $r_2$ denote the double roots of the Wronskian, the two
solutions of $T=0$ correspond to the following subspaces:
\begin{gather*}
  U = \lspan \{ (z-r_1) (z-r_2)^2, (z-r_1)^2 (z-r_2) \},\\
  U = \lspan \{ (z-r_1)^3, (z-r_2)^3 \}.
\end{gather*}
The f\/irst of these is projectively equivalent to $\lspan\{z,z^2\}$,
while the second is equivalent to $\lspan\{ 1,z^3\}$.

The analysis of case $(iii)$ requires us to consider the solutions of
the equation
\[ (16/9) H^3 + T^2 = 0.\]
Calculation with a Gr\"obner basis package yields the following
necessary condition:
\[ a_{11} + 3 a_{20} = 0.\]
Further computation reveals that the general solution is given by
\begin{gather*}
   a_{11}  = -b^3/4+(3/4) b^2 a_{21} - a_{21}^3,\qquad
  a_{20} = -b^2/4+a_{21}^2/4,\\
   U = \lspan\{(z+(a_{21}+3b)/2)^2(z+(a_{21}-b)/2),
  (z+(a_{21}-b)/2)^3 \},
\end{gather*}
where $b$ is an arbitrary, non-zero complex parameter.  Evidently, the
above subspace is projectively equivalent to $\{ 1, z^2\}$.

As our f\/inal example, let us consider a polynomial subspace  with a
f\/inite symmetry group.  More specif\/ically, let us consider subspaces with
Wronskian $W=z^4-1$.  Such a subspace must have a basis of the form
\[ U = \lspan \{ (z-1)^2 (z+a), (z+1)^2 (z+b)\}.\]
As expected, there are exactly two solutions:
\begin{gather*}
   a= e^{\pi i/3},\qquad b= e^{-2\pi i/3},\qquad U_1 = \lspan\{ z^3 -
  \sqrt{3} i z, z^2-i/\sqrt{3} \},\\
   a= e^{-\pi i/3},\qquad b= e^{2\pi i/3},\qquad U_2 = \lspan\{ z^3 +
  \sqrt{3} i z, z^2+i/\sqrt{3} \}.
\end{gather*}
Let us consider the question of whether~$U_1$ is projectively
equivalent to~$U_2$.  By construction, the two subspaces have the same
Wronskian.  We must therefore consider the projective
symmetries of $W=z^4-1$.  Applying the necessary formulas, we arrive at
\begin{gather*}
  J^2(z)  = -\frac{4(1+z^4)^2}{9 z^4},\qquad
  K(z)  = - \frac{2(1+10z^4+z^8)}{9 z^4}.
\end{gather*}
Both of the rational function $J^2(z)- J^2(Z)$ and $K(z) - K(Z)$ have the same
numerator, namely
\[ (z - Z)(z - iZ)(z + iZ)(z + Z)(-1 + zZ)(-i + zZ)(i + zZ)(1 + zZ). \]
Each of these factors def\/ines a fractional linear transformation.  The
union of all such forms the octahedral symmetry group of the
polynomial $z^4-1$.  By Theorem~\ref{thm:genequiv} the subspaces~$U_1$
and~$U_2$ are projectively equivalent if and only if they are related
by one of these transformations.  Indeed, exactly half of these
transformations, the subgroup consisting of
\[ z=Z, \qquad z=1/Z,\qquad z=-Z, \qquad z=-1/Z \]
are symmetries of $U_1$ and $U_2$.  The other coset, consisting of the
four transformations
\[ z=i Z,\qquad z=i/Z,\qquad z=-i Z,\qquad z=-i/Z \] relate $U_1$, $U_2$.  The two
subspaces in question are projectively equivalent.

Remark: it is an interesting question to characterize and if possible
to classify those polynomial subspaces whose projective symmetry group is
strictly smaller than the symmetry group of the corresponding Wronskian.

\subsection{The problem of top powers}
Let us say that a $p\in \cP_n$ is a top power if it is of the form
$p(z) =(az+b)^n$ where $(a,b)\neq (0,0)$.  The following result is
well-known in classical invariant theory~\cite{kungrota}.

\begin{theorem}[Sylvester]
  \label{thm:sylvester}
  Generically, a complex polynomial $($binary form$)$ having odd deg\-ree
  $n=2m+1$ can be given as a sum of~$m$ distinct top powers.
\end{theorem}

Indeed, determining the smallest number of $n$th powers required for
the expression a given $p\in \cP_n$ is a central problem of classical
invariant theory.  In this last subsection, we modify this question to
the context of polynomial subspaces.  First, we observe the following.

\begin{Proposition}
  Let $U\in\cG_k\cP_n$, where $k\leq n$, be a polynomial subspace.
  Then, $U$ contains at most $k$ top powers.
\end{Proposition}

\begin{proof}
  Suppose that $U$ contains $k+1$ top powers.  By applying an $\SL_2$
  transformation, if necessary, without loss of generality, none of
  these is constant. By Proposition~\ref{prop:npower1}, every element
  of $U^*$ is divisible by $q(z) = \prod_{i=1}^{k+1} (z-b_i)$ where
  the $b_i$ are distinct.  However, the subspace of $\cP_n$ of
  polynomials divisible by $q(z)$ is $(n-k)$-dimensional.  However,
  $\dim U^*= \ell = n+1-k$~-- a contradiction.
\end{proof}

Also observe that if $U\in \cG_k\cP_n$ contains $k$ distinct top
powers, then these act as a basis for~$U$.  Geometrically, a subspace
that admits a basis of top powers can be characterized as a secant
f\/lat of the rational normal curve~\cite{chipger}.  Here, we focus on
the following.

\begin{Question}
  How many top powers $($up to scalar multiple$)$ does a given polynomial
  subspace contain?  Formulate a procedure to decide whether a
  polynomial subspace admits a basis of top powers.
\end{Question}

The following procedure is one possible answer to this question.
For a polynomial subspace~$U$, let $U^{(j)}$, $j\geq 0$ denote the
subspace spanned by the $j$th derivatives of the elements of~$U$.

\begin{theorem}
  \label{thm:npower}
  Let $U$ be a $k$-dimensional polynomial subspace with
  $n=\deg U$ and $\ell=n+1-k$ the codimension.  Let $Q=
  \W(U^{(\ell-1)})$.  If $1\notin U$, then $U$ admits a basis of top
  powers if and only if the following conditions hold:
  \begin{itemize}\itemsep=0pt
  \item[$(i)$] $Q$ has $k$ simple roots   and
  \item[$(ii)$] every $q\in U^*$ is divisible by $Q$.
  \end{itemize}
\end{theorem}

\begin{proof}
  Suppose that $U$ contains $k$ distinct top powers. Since $1\notin U$
  these have the form $(z-b_i)^n$, $i=1,\ldots, k$ for some $b_i\in
  \Cset$.  Observe that
  \[Q(b_i)=\W((z-b_1)^{n-\ell+1},\ldots,
  (z-b_k)^{n-\ell+1})\big|_{z=b_i} = 0,\qquad i=1,\ldots, k\] and that
  $\deg Q=k$.  Condition (i) is established.  Condition $(ii)$ follows
  by Proposition~\ref{prop:npower1}.  Conversely, if conditions $(i)$
  and (ii) hold then by Proposition~\ref{prop:npower1}, $U$ contains
  $k$ distinct top powers.
\end{proof}

Mutatis mutandi we can formulate a version of the above
theorem for the case where $1\in U$.  One merely has to replace $k$
with $k-1$ in condition $(i)$.  Criterion $(i)$ can, in principle, be
tested by calculating the discriminant of~$Q$.  Criterion~$(ii)$
requires an explicit basis of $U^*$.  Such a~basis can be obtained by
means of Theorem~\ref{thm:dual1}.  Therefore, the above theorem can
serve as an algorithm for deciding the question of a top power basis.
We also note that one direction of Theorem~\ref{thm:npower}
follows from the more general Lemma~3.14 of~\cite{chipger}.

As an interesting special case, we have the following.   Let $U \in
\cG_k \cP_k$ be a codimension 1 subspace and $\W(U)\in \cG_1\cP_k$ the
corresponding Wronskian covariant.  Counting $\infty$, let
 $1\leq m\leq k$ be the number of
 distinct roots of $\W(U)$.
\begin{Proposition}
  \label{prop:codim1}
  The codimension $1$ subspace $U$ contains  $m$ top
  powers. In particular,  $U$
  admits a basis of top powers if and only if $\W(U)$ has simple roots
  only and $k-1\leq \deg \W(U)\leq k$.
\end{Proposition}
\begin{proof}
  By Theorem \ref{thm:Wdual}, $U^*$ is 1-dimensional, spanned by
  $\W(U)$.  The desired conclusion now follows by Propositions
  \ref{prop:npower0} and \ref{prop:npower1}.  In the generic case
  where $\W(U)\in \cP_k$ has $k$ distinct roots, $U$ admits a basis of
  the form $(z-b_i)^k$ where $b_1,\ldots, b_k$ are the distinct roots
  of $\W(U)$.  If $\deg W(U)=k-1$; i.e.\ if $\infty$ is a root of
  $\W(U)$, then the basis of top powers consists of $1$ and
  $(z-b_i)^k$ where $b_1,\ldots, b_{k-1}$ are the distinct roots of
  $\W(U)$.
\end{proof}

Our f\/inal contribution is the following conjecture.
\begin{Proposition}
    \label{prop:lthpower}
  If $U\in \cG_k\cP_n$ admits a basis of top powers then $W(U)$ is an
  $\ell$th power with
  \begin{equation}
    \label{eq:lthpower}
     W(U) = W\big(U^{(\ell-1)}\big)^\ell,
   \end{equation}
   where the equality is relative to $\cG_1 \cP_{k\ell}$; i.e.\ up to
   non-zero scalar multiple.
\end{Proposition}

\begin{Conjecture}
  Suppose that \eqref{eq:lthpower} holds and that $W(U^{(\ell-1)})$
  has $k$ distinct finite, simple roots.  Then, $U$ admits a basis of
  top powers.
\end{Conjecture}

Proposition \ref{prop:codim1} establishes the conjecture in the case
where the codimension $\ell=1$. We can also prove the conjecture for
the case of codimension $\ell=2$.  Since all of the roots of
$W(U^{(1)})$ are assumed to be f\/inite, $1\notin U$ and hence $U^{(1)}$
is $k$-dimensional. Let $b$ be one of the $k$ distinct roots of
$W(U^{(1)})$.  From the assumption that $W(U) = W(U^{(1)})^2$, it
follows that $b$ is a double root of~$W(U)$.  Let $\blambda$ be the
$b$-shape of $U$ and $\bnu=\blambda^+$ the corresponding pivot
sequence.  If $\nu_k=n$, we are done; for this means that $(z-b)^n\in
U$.  Suppose then that $\nu_k<n$.  Since $n=k+1$, this means that
either $\nu_k=n-2=k-1$ or that $\nu_k=n-1=k$.  If the former where
true then, $\nu_j = j-1$, $j=1,\ldots, k$ and hence $b$ is not a root
of~$W(U)$~-- a contradiction.  Let's consider the case where
$\nu_k=k$.  Hence, $U$ contains a polynomial of the form $(z-b)^k+ a
(z-b)^{k+1}$. Hence $U^{(1)}$ contains a polynomial of the form
$(z-b)^{k-1}+ a_1 (z-b)^{k}$.  Let $\tnu_i$ be the $b$-pivot sequence
of $U^{(1)}$.  Since $b$ is assumed to be a simple root of
$W(U^{(1)})$, we must have $\tnu_j = j-1$ for $j=1,\ldots, k-1$ and
$\tnu_k = k$.  However, we just argued that $U^{(1)}$ has a $b$-pivot
in column $k-1$; this is a contradiction.  Our claim is established.

\subsection*{Acknowledgements}
Discussions with D.~G\'omez-Ullate and N.~Kamran are
  gratefully acknowledged.  P.C.~was supported by an NSERC
  Undergraduate Summer Research Award. R.M.~is supported by an NSERC
  discovery grant.

\pdfbookmark[1]{References}{ref}
\LastPageEnding

\end{document}